%% file: bl-pow-sing.tex
%

\input def1200.tex

\overfullrule=0pt
{\nopagenumbers
\null
\vskip 4 true cm
\title{BLOWING UP THE POWER OF A SINGULAR CARDINAL}
\sect{Moti Gitik}
\ce{School of Mathematics}
\ce{Raymond and Beverly Sackler}
\ce{Faculty of Exact Sciences}
\ce{Tel Aviv University}
\ce{Tel Aviv 69978 Israel}
\vfill\eject}
\count0=1
\null
\dspace
\def\llvdash{\mathop{\|\hskip-2pt \raise 2pt\hbox{\vrule
height 0.15pt width 0.25cm}}}
\sect{Introduction}

Suppose that $\kap$  is a singular cardinal of
cofinality $\ome$.  We like to blow up its
power.  Overlapping extenders where used for
this purpose in [Git-Mag2].  On the other hand,
it is shown in [Git-Mit] that it is necessary to
have for every $n<\ome$ unboundedly many
$\alp$'s in $\kap$  with $o(\alp)\ge\alp^{+n}$.
The aim of the present paper is show that this
assumption is also sufficient.  Ideas of [Git
hid. et] will be extended in order to produce
$\kap^{++}$  $\ome$-sequences.  In [Git hid.
ext] an $\ome$-sequence corresponding to two
different sequences of measures was constructed.
Here we would like to construct a lot of
$\ome$-sequences corresponding to the same sequence 
of measures.

The first stage will be to to force with a
forcing which produces $\kap^{++}$  Prikry
sequences but the cost is that $\kap^{++}$  and
is collapsed.  Then a projection of this forcing
will be defined such that the resulting
forcing  will still have $\kap^{++}$  Prikry
sequences but also satisfy $\kap^{++}$-c.c. and
preserve $\kap$  strong limit cardinal.

\sect{1.~~Preparation Forcing-- the first try }

Let us assume GCH.  Suppose that $\kap_\ome
=\bigcup_{n<\ome}\kap_n$  and
$o(\kap_n)= \kap_n^{+n+2}+1$.
We will define a forcing which will combine
ideas of [Git-Mag2] and [Git hid. ext].  In
contrast to [Git hid. ext] we like to produce
lots of Prikry sequences even by the cost of
collapsing cardinals.  The main future of this
forcing will be the Prikry condition.  Splitting
it above and below $\kap_n$ $(n<\ome)$  we will
be able to conclude that the part above $\kap_n$
does not add new subsets to $\kap_n$  and the
part below does not effect cardinals above $\kap_n$.
The problematic cardinal will be $\kap_\ome^{++}$.
In order to prevent it from collapsing we
construct a projection of the forcing which will
satisfy $\kap^{++}_\ome$-c.c.

For every $n<\ome$.  Let us fix a nice system
$\UU_n=\ll \calU_{n,\alp}\mid\alp <\kap_n^{+n+2}>,
<\pi_{n,\alp,\bet}\mid\alp,\bet <\kap_n^{+n+2},\
\calU_{n,\alp}\triangleleft\calU_{n,\bet}\gg$.
We refer to [Git-Mag1] for the basic definitions.
Actually an extender of the length
$\kap_n^{+n+2}$ will be fine for our purpose as
well.

For every $n<\ome$, let us first define a
forcing notion $\langle Q_n,\le_n\rangle$  and
then use it as the level $n$  in the main
forcing.

Fix $n<\ome$.  We like to define a forcing
$\langle Q_n,\le_n\rangle$.   Let us drop the lower
index $n$ for a while.

$Q$  will be the union of two sets $Q^0$  and
$Q^1$  defined below.

\subheading{Definition 1.1}  Set $Q^1$  to be
the product of $\{ p\mid p$ is a partial
function from $\kap^{+n+2}$  to $\kap^{+n+2}$
such that $\dom p$  is an ordinal less than
$\kap^{+n+2}\}$ and $\{$  $q\mid q$  is a
partial function from $\kap_\ome^{++}$  to
$\kap^{+n+2}$  of cardinality less than
$\kap^+_\ome\}$.

The ordering on $Q^1$ is an inclusion.  I.e.
$Q^1$ is the product of the product of two Cohen
forcings:  for adding a new subset to $\kap^{+n+2}$
and for adding $\kap_\ome^{++}$  new subsets to
$\kap_\ome^+$.

\subheading{Definition 1.2}  A set $Q^0$
consists of triples $\langle p,a,f\rangle$
where
\item{(1)} $p=\langle\{ <\gam, p^\gam>\mid\gam <\del\},
g,T\rangle$ where
\itemitem{(1a)} $g\subseteq\kap^{+n+2}$  of cardinality
$<\kap$.
\itemitem{(1b)} $\del<\kap^{+n+2}$ 
\itemitem{(1c)} $o\in g$  and every initial segment of $g$
(including $g$  itself) has the least upper
bound in $g$.
\itemitem{(1d)} $\del > \max (g)$	
\itemitem{(1e)} for every $\gam\in g$ $p^\gam$  is the
empty sequence 
\itemitem{(1f)} $T\in\calU_{\max (g)}$
\itemitem{(1h)} for every $\gam\in\del\bks g$  $p^\gam$
 is an ordinal below $\kap_\ome^{++}$.
\itemitem{}
\item{}

Further we shall denote $g$  by $\supp (p)$, the
maximal element of $g$  by $mc(p)$,  $\del$  by
$\del (p)$ and $T$  by $T(p)$.  Let us refer to
ordinals below $\del(p)$  as coordinates. We
will frequently confuse between an ordinal
$\gam$ and one element sequence $\langle
\gam\rangle$.  

\item{(2)} $a$  is a partial one to one order
preserving function between $\kap^{++}_\ome$
and $\del(p)$  of cardinality less than $\kap$.
Also every $\gam\in\dom a$  is below $mc(p)$  in
sense of the ordering of extender $\UU$.
\item{(3)} $f$  is a partial function from
$\kap_\ome^{++}$  to $\kap^{+n+2}$ of cardinality
less than $\kap_\ome^+$  and such that $\dom
f\cap \dom a=\emptyset$.
\item{}

Let us give some intuitive motivation for the
definition of $Q^0$.  Basically we like to add
$\kap_\ome^{++}$.  Prikry sequences (actually
a one element sequence).

The length of the extender used is only
$\kap^{+n+2}$.  A typical element of $Q^0$
consists of a triple $\langle p,a,f\rangle$.
The first part of it $p$  is as a condition of
[Git-Mag1] with slight changes need for mainly
technical reasons.  The idea is to assign
ordinals $<\kap_\ome^{++}$  to the coordinates
of such $p$'s.  $a$ is responsible for this
assignment.  Basically, if for some
$\alp<\kap_\ome^{++}$, $\bet <\kap^{+n+2}$
$a(\alp)=\bet$,  then $\alp$-th sequence will be
read from the $\bet$-th Prikry sequence.
Clearly, we do not want to allow this assignment
to grow into the one to one correspondence between
$\kap^{+n+2}$  and $\kap_\ome^{++}$.  The third
part $f$  and mainly the definition of the
ordering below is designed to prevent such
correspondence.  

\subheading{Definition 1.3}  $Q=Q^0\cup Q^1$.

Let us turn to the definition of the order over
$Q$.  First we define $\le^*$ the pure extension.

\subheading{Definition 1.4} Let $t,s\in Q$.
Then $t\le^*s$ if either
\item{(1)} $t,s\in Q^1$  and $t$  is weaker than $s$
in the ordering of $Q^1$  or
\item{(2)} $t,s\in Q^0$  and the following holds:\hb
let $t=\langle p,a,f\rangle$, $s=\langle
q,b,g\rangle$ (2a) $p\le^* q$  in sense of [G2t-Mag1]
with only addition in (v):
\item{(i)} $\del (p)\le \del (q)$
\item{(ii)} $\supp(p) \subseteq \supp(q)$ 
\item{(iii)} for every $\gam<\del(p)$
$p^\gam=q^\gam$
\item{(iv)} $\pi_{mc(q)mc(p)}$  projects $T(q)$
into $T(p)$
\item{(v)} for every $\gam\in\supp (p)\cup \dom
a$  and $\nu\in T(q)$ 
$$\pi_{mc(q),\gam}(\nu)=\pi_{mc(p),\gam}(\pi_{mc(q),
mc(p)}(\nu))\ .$$
\itemitem{(2b)} $a\subseteq b$
\itemitem{(2c)} $f\subseteq g$.
\item{}

Notice that in contrast to [Git-Mag1], the
commutativity in (2a)(v) does not cause a
special problem since the number of coordinates
$\supp (p)\cup\dom a$  has cardinality $<\kap$,
i.e. below the degree of completeness of ultrafilters
in the extender used here.  

\subheading{Definition 1.4.1} Let $s,t\in Q$.  We
say that $s$  extends $t$  if $t\le^*s$  or
$t\in Q^0$, $s\in Q^1$  and the conditions below
following hold.

Let $t=\langle p,a,f\rangle$ and $s=\langle
q,h\rangle$. 
\item{(1)} $\del(p)\le\dom q$  (recall that by 1.1,
$\dom q$  is an ordinal $<\kap^{+n+2}$).
\item{(2)} for every $\gam\in\del (p)\bks \supp (p)$
if $p^\gam <\kap^{+n+2}$  then $p^\gam =q(\gam)$
otherwise $q(\gam)=\kap$.
\item{(3)} $q(mc(p))\in T(p)$
\item{(4)} for every $\gam\in\supp (p)$ $q(\gam)=
\pi_{mc(p),\gam} (q(mc(p)))$
\item{(5)} $h\supseteq f$
\item{(6)} $\dom h\supseteq \dom a$ 
\item{(7)} for every $\bet \in \dom a$ 
$h(\bet)=q (a(\bet))$,  if $a(\bet)\in\supp (p)$
or $h(\bet)=\pi_{mc(p), a(\bet)}(q(mc(p)))$,  otherwise.
\item{}

The conditions (1) to (4) are as in [Git-Mag 1]
with only change in (2) in case $p^\gam\ge \kap^{+n+2}$.
Then it is replaced by $\kap$.  The idea behind
this is to remove unnecessary information a 
condition may have in order to prevent collapses
of cardinals above $\kap^{+n+2}$.  The
conditions (5) to (7) are the heard of the
matter.  Our purpose is to forbid the assignment
$a$ from growing into a $1-1$ function from 
$\kap_\ome^{++}$  to $\kap^{+n+2}$  but to still  
produce $\kap_\ome^{++}$-sequences.  What
actually happens in the definition is a switch
from Prikry type harmful forcing to a nice Cohen
type forcing.  The only essential information
from $a$  is put into $h$.  The actual place of
the sequence $\bet (\bet\in\dom a)$  is hidden 
after passing from $t$  to $s$.   

\proclaim Lemma 1.5.  $Q^1$  is dense in $Q$.

The proof follows from  Definition 1.4.1.

\proclaim Lemma 1.6.  $\langle Q,\le\rangle$
does not collapse cardinals or blows up their
powers. 

Follows from 1.5.

\proclaim Lemma 1.7.  $\langle Q,\le,\le^*\rangle$
satisfies the Prikry condition.

The proof of the parallel statement of [Git-Mag 1]
applies here without essential changes.

Now let us put all $Q_n$'s defined above
together.

\subheading{Definition 1.8}  A set of forcing
conditions $\calP$  consists of all elements $p$
of the form $\langle p_n\mid n<\ome\rangle$  so
that
\item{(1)} for every $n<\ome$ $p_n\in Q_n$
\item{(2)} there exists $\ell <\ome$  such that for
every $n\ge\ell$  $p_n\in Q^0_n$.

Let us denote further the least such $\ell$  by
$\ell (p)$.

\subheading{Definition 1.9}  Let $p=\langle
p_n\mid n<\ome\rangle$, $q=\langle q_n\mid n<\ome
\rangle\in\calP$.  We say that $p$  extends
$q(p\ge q)$  if for every $n<\ome$  $p_n$
extends $q_n$  in the ordering of $Q_n$.  

\subheading{Definition 1.10}  Let $p,q\in\calP$.
We say that $p$  is a direct or pure extension
$q$ iff $p\ge q$  and $\ell (p)=\ell (q)$.

\proclaim Lemma 1.11.  $\langle
\calP,\le,\le^*\rangle$  satisfies the Prikry
condition.

\subheading{Sketch of the Proof}
Let $\sig$  be a statement of the forcing
language and $p\in\calP$.  We are looking for
$q\ge^*p$  deciding $\sig$.  Assume for
simplicity that $\ell (p)=0$.  As in [Git-Mag 1]
we extend $p$  level by level trying to decide
$\sig$.  Suppose that we passed  level 0  and
are now on level 1.  We have here basically two
new points.  The first to our advantage is that
the measures on the level 1 are $\kap_1$-complete
and $\kap_1>\kap_0$.  So we can always shrink
sets of measure 1 in order to have the same
condition in $Q^0_0$ on the level 0.  The second
point is that the cardinality of $Q^1_0$  is
big. However let us then use the completeness of
$Q^1_0$.  Recall that $Q^1_0$  is
$\kap_\ome^+$-closed forcing.

The rest of the proof is parallel to [Git-Mag
1].\hfill$\square$

Let $G$  be a generic subset of $\calP$. For
$\bet <\kap^{++}_\ome$  let $G(\bet): \ome
\to\kap_\ome$ be the function defined as
follows.  $G(\bet)(n)=\nu$  iff there is
$\langle p_k\mid k<\ome\rangle\in G$  such that
$\bet\in\dom p_{n,2}$  $p_{n2}(\bet)=\nu$,
where $p_{n,2}$  is the second coordinate of
$p_n\in Q^1_n$.   

Notice that we cannot claim $G(\bet)$'s are
increasing with $\bet$.  Actually, lots of them
will be old sequences and also they may be equal
or reverse the order.  But the following is still
true.

\proclaim Lemma 1.12.  For every $\gam <\kap_\ome^{++}$ 
there is $\bet$, $\gam <\bet <\kap_\ome^{++}$
such that $G(\bet)$  is above every $G(\bet')$
with $\bet'<\bet$.

\pr Work in $V$.  Let $p\in\calP$.  Suppose for
simplicity that $\ell (p)=0$. Otherwise work
above the level $\ell (p)-1$.  Let $p=\langle
p_n\mid n<\ome\rangle$  and $p_n=\langle
p_{n0},p_{n1}, p_{n2}\rangle$  $(n<\ome)$.
Pick some $\bet$,  $\gam < \bet <\kap_\ome^{++}$
which above everything appears in $p$, i.e.
$\bet >\cup \{ \del (p_{n0})\cup\sup (\dom
p_{n1}\cup\dom p_{n2})\mid n<\ome\}$.
Extend $p$  to a condition $q=\langle q_n\mid
n<\ome\rangle$, $q_n=\langle q_{n0}, q_{n1},q_{n2}
\rangle$ such that $q_{n1}=p_{n1},
q_{n2}=p_{n2}$ and $mc(q_{n0})>mc(p_{n0})$  for
every $n<\ome$.  Extend now $q$  to $r=\langle
r_n\mid n<\ome\rangle$,  $r_n=\langle
r_{n0},r_{n1},r_{n2}\rangle$ by adding the pair
$\langle \bet, mc(q_{n0})\rangle$  to $q_{n1}$
for every $n<\ome$. 

We claim that
$$r\llvdash\left(\buildrul\sim\under
G(\bet)>\buildrul\sim \under G(\bet')\quad\hbox{for every}
\quad\bet'<\bet)\right)\ .$$ 

Fix $\bet'<\bet$  and let $s\ge r$.  W.l. of $g$.
$\ell (s)=\ell (r)=0$.  Since otherwise we
repeat the same argument above $\ell (s)$.  Let
$s=\langle s_n\mid n<\ome\rangle$ and
$s_n=\langle s_{n0},s_{n1},s_{n2}\rangle$ for
every $n<\ome$.  Denote by $A$  the set of all
$n$'s such that $\bet'\in\dom s_{n1}$.  For
every $n\in\ome\bks A$  extend $s_n$ by adding
there pair $\langle\bet',0\rangle$  to $s_{n2}$.
Let us still denote the resulting condition by
$s$.  Then the function $G(\bet')\rhookup
\ome\bks A$  will be forced by $s$  to be an old
function.  Hence $G(\bet)\rhookup \ome\bks A$
is above it.   

Now let $n\in A$.  Then, since $\bet'<\bet$,
$\bet',\bet\in\dom s_{n1}$  and $s_{n1}$  is
order preserving, the coordinate assigned to 
$\bet'$  by $s_{n1}$  is below the one assigned
to $\bet$.  Hence $s$  forces that $\buildrul\sim\under
G(\bet)\rhookup A$ is above $\buildrul\sim\under
G (\bet')\rhookup A$  and we are done.\hfill$\square$ 

For $n<\ome$ let us split $\calP$  into
$\calP\rhookup n$  and $\calP\bks n$  as follows: 
$$\eqalign{&\calP\rhookup n=\{ p\rhookup n\mid
p\in\calP\}\cr
&\calP\bks n=\{ p\bks n\mid p\in\calP\}\
.\cr}$$
The following lemma is routine

\proclaim Lemma 1.13.  For every $n<\ome$  the
forcing with $\calP$  is the same as the forcing
with $(\calP\bks n)\times(\calP\rhookup n)$.  

\proclaim Lemma 1.14.  $\langle \calP,\le\rangle$
preserves the cardinals $\le\kap^+_\ome$  and
GCH holds below $\kap_\ome$  in a generic
extension by $\calP$.

\pr  For every $n<\ome$ $\kap_{n+1}$  is
preserved since $\calP$  splits as 1.13  into a
forcing $\calP\bks n$  and $\calP\rhookup n$.
By analogous of 1.11 for $\calP\bks n$,
$\calP\bks n$  does add new bounded subsets of
$\kap_{n+1}$.  By 1.6, $\calP\rhookup n$
preserves cardinals.   
Therefore, nothing below $\kap_\ome$ is
collapsed.  Now if $\kap_\ome^+$  is collapsed
then $|\kap^+_\ome|=\kap_\ome$  which is
impossible by the Weak Covering Lemma
[Mit-St-Sch] or just directly using arguments
like those of [Git-Mag 1], Lemma
1.11.\hfill$\square$   

Unfortunately, $\kap_\ome^{++}$  is collapsed
by $\calP$  as it is shown in the next lemma.

\proclaim Lemma 1.15.  In $V[G]$  $\ |(\kap_\ome^{++})^V
|=\kap^+_\ome$.

\pr Work in $V$.  The cardinality of the set
$\prod_{n<\ome}\kap_n^{+n+2}/$  finite is
$\kap^+_\ome$.  Fix some enumeration $\langle
g_i\mid i<\kap_\ome^+\rangle$  of it. 

Now in $V[G]$,  let $p=\langle p_n\mid
n<\ome\rangle$  $\in G$, $p_n=\langle
p_{n0},p_{n1},p_{n2}\rangle$  $(n<\ome)$,  $\bet
<\kap_\ome^{++}$  and starting with some
$n_0<\ome$  $\bet\in\dom p_{n1}$.  Find $i<\kap^+_\ome$	
s.t. the function $\{\langle n,
p_{n1}(\bet)\rangle\mid n\ge n_0\}$ belongs to
the equivalence class $g_i$. Set then
$i\mapsto \bet$.  Using genericity of $G$  it is
easy to see that this defines a function from
$\kap^+_\ome$  unboundedly into $\kap_\ome^{++}$.
\hfill$\square$ 
	
We would like to project the forcing $\calP$  to
a forcing preserving $\kap_\ome^{++}$.  The idea
is to make it impossible to read from the
sequence $G(\bet)$ $(\bet <\kap_\ome^{++})$  the
sequence of coordinates ($\MOD$  finite) which
produces $G(\bet)$  in sense of 1.15.  The
methods of [Git] will be used for this
purpose.  But first the forcing $\calP$ should be
fixed slightly.  The point is that we like to
have much freedom in moving $\bet$'s from the beginning.
$\calP$ is quite rigid in this sense. Thus, for
example, if some $\bet <\kap_\ome^{++}$
corresponds to a sequence of coordinates $g$  in
$\prod_{n<\ome}\kap^+_n$,  then using $G(\bet)$
only it is easy to reconstruct $g$  modulo finite.

\sect{2.~~The Preparation Forcing}

Suppose that $n<\ome$  is fixed.  For every
$k\le n$  we consider a language $\calL_{n,k}$
containing a constant $c_\alp$ for every $\alp
<\kap_n^{+k}$ and a structure 
$$\gra_{n,k}=\langle H(\lam^{+k}), \in, \lam,
0,1\nek \alp\nek\mid\alp <\kap_n^{+k}\rangle$$ 
in this language, where $\lam$  is a regular
cardinal big enough.  For an ordinal $\xi
<\lam$  (usually $\xi$  will be below
$\kap_n^{+n+2})$  we denote by $tp_{n,k}(\xi)$
the $\calL_{n,k}$-type realized by $\xi$  in
$\gra_{n,k}$.  Let $\del <\lam$.  $\calL_{n,k,\del}$
will be the language obtained from $\calL_{n,k}$
by adding a new constant $c$.  $\gra_{n,k,\del}$
will be $\calL_{n,k,\del}$-structure obtained
from $\gra_{n,k}$  by interpreting $c$  as $\del$.
The type $tp_{n,k}(\del,\xi)$ is defined in the
obvious fashion.  Further we shall freely identify
types with ordinals corresponding to them in
some fixed well ordering of the power sets of
$\kap_n^{+k}$'s.  The following is an easy
statement proved in [Git].

\proclaim Lemma 2.0.  Suppose that
$\alp_0,\alp_1 <\kap_n^{+n+2}$  are realizing
the same $\calL_{n,k,\rho}$-type for some $\rho
<\min (\alp_0,\alp_1)$ and $n\ge k>0$.  Then for
every $\bet$,  $\alp_0\le\bet <\kap_n^{+n+2}$
there is $\gam, \alp_1\le \gam <\kap_n^{+n+2}$
such that the $k-1$-type realized by $\bet$ over
$\alp_0$ (i.e. $\calL_{n,k-1,\alp_0}$-type) is
the same as those realized by $\gam$  over $\alp_1$.

\proclaim Lemma 2.1.  Let $\gam <\kap_n^{+n+2}$.
Then there is $\alp <\kap_n^{+n+2}$  such that
for every $\bet\in (\alp,\kap_n^{+n+2})$  the
type $tp_{n,n}$  $(\gam,\bet)$  appears (is
realized) unboundedly often in $\kap_n^{+n+2}$.

\pr  The total number of such types is
$\kap_n^{+n+1}$.  Let $\langle t_i\mid
i<\kap_n^{+n+1}\rangle$  be an enumeration of
all of them.  For each $i<\kap_n^{+n+1}$  set
$A_i$  to be the subset of $\kap_n^{+n+2}$
consisting of all the ordinals realizing $t_i$.
Define $\alp$ to be the supremum of $\{\cup
A_i\mid i<\kap_n^{+n+1}$ and $A_i$ is bounded in
$\kap_n^{+n+2}\}$.\hfill$\square$ 

\proclaim Lemma 2.2.  Let $\gam <\kap_n^{+n+2}$.
Then there is a club $C\subseteq\kap_n^{+n+2}$
such that for every $\bet\in C$  the type
$tp_{n,n}(\gam,\bet)$  is realized stationary
many times in $\kap_n^{+n+2}$. 

\pr Similar to 2.1.\hfill$\square$

\proclaim Lemma 2.3. The set $C=\{\bet
<\kap_n^{+n+2}|$  for every $\gam <\bet$ 
$\ tp_{n,n}(\gam,\bet)$  is realized stationary
often in $\kap_n^{+n+2}\}$ containing a club.

\pr  Suppose otherwise. Let $S=\kap_n^{+n+2}\bks C$.
Then
$$S=\{\bet <\kap_n^{+n+2}\mid \exists \gam
<\bet\ tp_{n,n}(\gam,\bet)\quad \hbox{appears only
nonstationary often in}\quad\kap_n^{+n+2}\}$$
and it is stationary.  Find $S'\subseteq S$ stationary and
$\gam' <\kap_n^{+n+2}$  such that for every
$\bet\in S'$  $tp_{n,n}(\gam',\bet)$  appears
only nonstationary often in $\kap_n^{+n+2}$.
But this contradicts 2.2. Contradiction.\hfill$\square$

For $\ell\le k\le n$ and $\calL_{n,k}$-type $t$
let us denote by $t\rhookup \ell$  the reduction
of $t$  to $\calL_{n,\ell}$,  i.e. the
$\calL_{n,\ell}$-type obtained from $t$  by
removing formulas not in $\calL_{n,\ell}$.

\proclaim Lemma 2.4.  Let $0 <k,\ell\le n$,
$\gam <\bet <\kap_n^{+n+2}$  and $t$  be a
$\calL_{n,\ell,\gam}$-type realized above $\gam$.
Suppose that $tp_{n,k}(\gam,\bet)$  is realized
unboundedly often in $\kap_n^{+n+2}$.  Then
there is $\del$, $\gam <\del<\bet$  realizing
$t\rhookup \min (k-1,\ell)$.

\pr Pick some $\alp$, $\gam <\alp<\kap_n^{+n+2}$
realizing $t$.  Let $\rho >\max (\bet,\alp)$  be
an ordinal realizing $tp_{n,k}(\gam,\bet)$.
Then $\rho$  satisfies in $H(\lam^{+k})$  the
following formula of $\calL_{n,k,\gam}$:
$$\exists y(c<y<x)\wedge
(H(\lam^{+k-1})\ \hbox{satisfies}\quad\psi (y)\
\hbox{for every}\ \psi\ \hbox{in the set of formulas
coded by}\ c_{t\rhookup\min (k-1,\ell)})\ .$$
Hence the same formula is satisfied by $\bet$.
Therefore, there is $\del$,  $\gam<\del <\bet$
realizing $t\rhookup \min (k-1,\ell)$.\hfill$\square$

The above lemma will be used for proving
$\kap_\ome^{++}$-c.c. of the final forcing via
$\Del$-system argument.

Let us specify now ordinals which will be allowed
further to produce Prikry sequences.

\subheading{Definition 2.5}  Let $k\le n$  and
$\bet <\kap_n^{+n+2}$.  $\bet$  is called
$k$-good iff\hb
(1) for every $\gam <\bet$
$tp_{n,k}(\gam,\bet)$  is realized unboundedly
many times in $\kap_n^{+n+2}$\hb
and
$$cf\bet\ge \kap_n^{++}\ .\leqno(2)$$

$\bet$  is called good iff for some $k\le n$  $\bet$
is $k$-good.

By Lemma 2.3, there are stationary many $n$-good
ordinals.  Also it is obvious that $k$-goodness
implies $\ell$-goodness for every $\ell\le k\le
n$.

\proclaim Lemma 2.5.1.  Suppose that $n\ge k>0$
and $\bet$  is $k$-good.  Then there are
arbitrarily large $k-1$-good ordinals below
$\bet$.

\pr Let $\gam <\bet$.  Pick some $\alp >\bet$
realizing $tp_{n,k}(\gam,\bet)$.  The fact that
$\gam <\bet<\alp$  and $\bet$ is $k-1$-good can  
be expressed in the language $\calL_{n,k,\gam}$
as in Lemma 2.4.  So they are in $tp_{n,k}(\gam,\bet)$.
Hence there is $\del$, $\gam <\del <\bet$  which
is $k-1$-good.\hfill$\square$

Let us now turn to fixing of the forcings
introduced in Section 1. 
We are going to use on the level $n$ a forcing
notion $Q^*_n$.  It is defined as $Q_n$  was with only
one addition that each ordinal in the range of assignment
functions is good.

\subheading{Definition 2.6}  A set $Q_n^*$  is
the subset of $Q_n$  consisting of $Q^1_n$  and
all the triples $\langle p,a,f\rangle$ of
$Q^0_n$  such that every $\alp\in rnga$ is
good.  The ordering of $Q^*_n$  is just the
restriction of the ordering of $Q_n$.

Lemma 1.5, 1.6 and 1.7 hold easily with $Q_n$
replaced by $Q^*_n$.  Let us show few additional
properties of $Q^*_n$  which are slightly more
involved.

\proclaim Lemma 2.7.  Suppose $\langle
p,a,f\rangle$ $\in Q^*_n$  and $\kap_\ome^{++}>\bet
>\sup (\dom a \cup \dom f)$.  Then there is a
condition $\langle q,b,f\rangle \ge^* \langle
p,a,f\rangle$  such that $\bet\in \dom b$  and
$b(\bet)$  is $n$-good.

\pr Using Lemma 2.3 find some $\xi
<\kap_n^{+n+2}$  above $mc(p)$  which is
$n$-good.  Now extend $p$  to $q$  such that
$\xi\in\supp(q)$.  Let $b=a\cup\{
\langle\bet,\xi\rangle\}$.  Then $\langle
q,b,f\rangle$  is as desired.\hfill$\square$

\proclaim Lemma 2.8.  Suppose that $\langle p,a,f\rangle$,
$\langle q,b,g\rangle \in Q^*_n$, $\bet\in\dom a$ it is   
$k$-good for $k>1$, $\{\gam_i\mid i<\mu\}\subseteq
(\bet\cap\dom b)\bks\dom f,\gam_0 >\sup (\bet\cap
\dom a)$  and $b(\gam_0)>\sup a\tagg$ $(\bet
\cap\dom a)$.  Then there is $\langle p^*,a^*,f\rangle$
a direct extension of $\langle p,a,f\rangle$
such that 
\item{(1)} $\{\gam_i\mid i<\mu\}\subseteq \dom
a^*$.
\item{(2)} for every $i<\mu$  $a^*(\gam_i)$ and
$b(\gam_i)$  are realizing the same $k-1$-type
\item{(3)} for every $i<\mu$,  if $b(\gam_i)$
is $\ell$-good $(\ell \le n)$  then
$a^*(\gam_i)$  is $\min (\ell, k-1)$-good. 
\item{(4)} if $t$  is the $n$-type over $\sup
(a\tagg(\bet \cap \dom a))$  realized by the
ordinal coding $\{ b(\gam_i)\mid i<\mu\}$,
then the code of $\{ a^*(\gam_i)\mid i<\mu\}$
realizes $t\rhookup k-1$.

\pr Denote $\sup (a\tagg (\bet\cap \dom a))$  by
$\rho$.  Let $t$  be the $n$-type over $\rho$
realized by the ordinal coding $\{
b(\gam_i)\mid i<\mu\}$.  By Lemma 2.4, there is
$\del$,  $\rho <\del <\bet$  realizing
$t\rhookup k-1$.  Let $\langle \xi i\mid i<\mu\rangle$
be the sequence coded by $\del$.  Define  
$$a^*=a\cup \{\langle \gam_i,\xi_i\rangle \mid
i<\mu\}\ ,\ p^*=p$$
and $f^*=f$.  Then $\langle p^*,a^*,f^*\rangle$
is as required.\hfill$\square$

\proclaim Lemma 2.8.1.  Suppose that $\langle
p,a,f\rangle$,  $\langle q, b,g\rangle \in
Q^*_n$  and $\bet\in\dom a,\ \gam\in\dom b$
are such that 
\item{(1)} $\bet$  is $k$-good for some $k\ge 2$
\item{(2)} $\bet\cap\dom a=\gam\cap\dom b$  and
for every $\del\in\bet\cap\dom a$
$a(\del)=b(\del)$
\item{(3)} $\bet >\sup (\dom b)$.
\item{}

Then there direct extensions $\langle
p^*,a^*,f\rangle\ge^* \langle p,a,f\rangle$ and
$\langle q^*,b^*,g\rangle$  $\ge^*$  $\langle q,b,g
\rangle$ such that 
\item{(a)}$\dom a^*=\dom b^*=\dom a\cup\dom b$
\item{(b)} for every $\del\in \dom a^*$
$a^*(\del)$  and $b^*(\del)$  are realizing the
same $k-2$-type over $\rho =_{df}\sup a\tagg
((\bet\cap \dom a))$
\item{(c)} for every $\del\in\dom b$  if $b(\del)$
is $\ell$-good then $a^*(\del)$  is $\min (\ell,
k-2)$-good   
\item{(d)} for every $\del\in\dom a$  if
$a(\del)$  is $\ell$-good then $b^*(\del)$  is
$\min (\ell, k-2)$-good
\item{(e)} $mc(p^*)$  and $mc(q^*)$  are
realizing the same $k-2$-type over $\rho$, more
over for every $\del\in \dom a\cup\dom b$  the
way $mc(p^*)$  projects to $a^*(\del)$  is the
same as $mc(q^*)$  projects to $b^*(\del)$.

\pr Let $s$  denotes the $k-1$-type realized by
$mc(q)$  over $\rho=\sup (a\tagg (\bet \cap \dom
a))$.  By Lemma 2.4, there is $\del$,  $\rho
<\del <\bet$  realizing $s$.  For every $\eta
\in \dom b$  let $\tileta$  be the ordinal
projecting from $\del$  exactly the same way as
$b(\eta)$  projects from $mc(q)$.  Notice that
for $\eta\in\dom b\cap\dom a$
$\tileta=b(\eta)=a(\eta)<\rho$.  Also,
$\tileta$  and $b(\eta)$  are realizing the same
$k-1$-type over and if $b(\eta)$ is $\ell$-good
then $\tileta$  is $\min(\ell, k-1)$-good, for
every $\eta\in\dom b$.

Pick $p^*$  to be a direct extension of $p$
with $mc(p^*)$  above $mc(p),\del$.  Set
$a^*=a\cup \{\langle \eta,\tileta\,\rangle
\mid\eta\in\dom b\}$.  Now we should define the
condition $\langle q^*,b^*,g\rangle$.  Since
$\del$  and $mc(q)$  are realizing the same $k-1$-type, by
Lemma 2.0 there exists $\nu$  realizing over
$mc(q)$  the same $k-2$-type as $mc(p^*)$ is
realizing over $\del$.  For $\eta\in\dom a$
define $\tileta$  as above only using $mc(p^*)$
and $\nu$  instead of $\del$  and $mc(q)$.   
Set $b^*=b \cup \{\langle \eta, \tileta\,\rangle
\mid\eta\in\dom a\}$.  Let $q^*$  be the
condition obtained from $q$  by adding $\nu$  as
a new maximal coordinate.  Then $\langle
q^*,b^*,g\rangle$  is as
desired.\hfill$\square$ 

Let us now define the forcing $\calP^*$.

\subheading{Definition 2.9} A set of forcing
conditions $\calP^*$ consists of all elements
$p=\langle p_n\mid n<\ome\rangle\in\calP$  such
that for every $n<\ome$
\item{(1)} $p_n\in Q_n^*$
\item{(2)} if $n\ge \ell (p)$  then $\dom
p_{n,1}\subseteq\dom p_{n+1,1}$ where
$p_n=\langle p_{n0},p_{n1},p_{n2}\rangle$
\item{(3)} if $n\ge \ell(p)$ and $\bet\in\dom
p_{n,1}$  then for some nondecreasing converging
to infinity sequence of natural numbers $\langle
k_m\mid \ome >m\ge n\rangle$  for every $m\ge n$
$p_{m,1}(\bet)$  is $k_m$-good.

The ordering of $\calP^*$  is as that of
$\calP$.

The intuitive meaning of (3) is that we are
trying to make the places assigned to the $\bet$-th
sequence more and more indistinguishable while
climbing to higher and higher levels.

The following lemma is crucial for transferring
the main properties of $\calP$  to $\calP^*$.

\proclaim Lemma 2.10.  $\langle
\calP^*,\le^*\rangle$  is $\kap_0$-closed.

\pr Let $\langle p(\alp)\mid \alp <\mu
<\kap_0\rangle$ be a $\le^*$-increasing sequence
of conditions of $\calP^*$.  Let for each
$\alp<\mu$  $p(\alp)=\langle p(\alp)_n\mid
n<\ome\rangle$  and for each $n<\ome$  $p(\alp)_n=\langle
p(\alp)_{n0}$, $p(\alp)_{n1}$, $p(\alp)_{n2}\rangle$.
For every $n<\ome$  find $q_{n0}\in Q_n^{0*}$
such that $q_{n0}\ge^*p(\alp)_{n0}$  for every
$\alp <\mu$.  Set $q_{n1}=\bigcup_{\alp<\mu}p(\alp)_{n,1}$
and $q_{n2}=\bigcupl_{\alp <\mu}p(\alp)_{n,2}$ for every
$n<\ome$.  Set $q_n=\langle q_{n0},q_{n1},q_{n2}\rangle$  
$(n<\ome)$ and $q=\langle q_n\mid n<\ome \rangle$.
Then $q\in\calP^*$.  Let us check the condition (3)
of Definition 2.9. Suppose that $\bet\in\dom
q_{n,1}$ for some $n<\ome$.  Then there is
$\alp<\mu$  such that $\bet\in \dom
p(\alp)_{n,1}$.  But now the sequence $\langle
k_m\mid \ome >m\ge n\rangle$  witnessing (3) for
$p(\alp)$  will be fine also for $q$.\hfill$\square$

Analogous of Lemmas 1.11, 1.13 and 1.14 hold for
$\calP^*$.  We define $\calP^*\rhookup n$ and
$]\calP^*\bks n$  from $\calP^*$  exactly
as $\calP\rhookup n$  and $\calP\bks n$  were
defined from $\calP$.

\proclaim Lemma 2.11. $\langle\calP^*,\le ,\le^*\rangle$
satisfies the Prikry condition.

\proclaim Lemma 2.12.  For every $n <\ome$ the
forcing with $\calP^*$  is the same as the
forcing with $(\calP^*\bks n)\times
(\calP^*\rhookup n)$.  

\proclaim Lemma 2.13.  $\langle\calP^*,\le\rangle$
preserves the cardinals below $\kap_\ome$  and
GCH below $\kap_\ome$  still holds in a generic
extension by $\calP^*$.

Let us show that $\calP^*$ adds lot of Prikry sequence.
Let $G$  be a generic subset of $\calP$.  For
$\bet <\kap_\ome^{++}$  we define
$G(\bet):\ome\to\kap_\ome$  as in Section 1,
i.e. $G(\bet)(n)=\nu$ iff there is $\langle
p_k\mid k<\ome >\in G$  such that $\bet\in\dom
p_{n,2}$  and $p_{n,2}(\bet)=\nu$  where
$p_n=\langle p_{n1},p_{n2}\rangle\in Q_n^{1*}$.

We claim that for unboundedly many $\bet$'s
$G(\bet)$  will be a Prikry sequence and
$G(\bet)$ will be bigger (modulo finite) than
$G(\bet')$  for every $\bet'<\bet$.  The next
lemma proves even slightly more.

\proclaim Lemma 2.14.  Suppose $p=\langle
p_k\mid k<\ome\rangle$  $\in\calP^*$,
$p_k=\langle p_{k0},p_{k1},p_{k2}\rangle$  for
$k\ge\ell (p)$, $\bet <\kap_\ome^{++}$  and
$\bet\notin \bigcup_{\ell(p)\le k<\ome}(\dom p_{k1}
\bigcup\dom p_{k2})$.  Then there is a direct extension
$q$  of $p$ such that $\bet\in\bigcup_{k\ge\ell (q)}\dom
q_{k,1}$, where $q=\langle q_k\mid
k<\ome\rangle$  and $q_k=\langle q_{k0},q_{k1},q_{k2}
\rangle$ for every $k\ge \ell (q)$.

\pr Let us assume for simplicity that $\ell(p)=0$.
Set $a=\bigcup_{k<\ome}\dom p_{k1}$.

\subheading{Case 1} $\bet\ge \bigcup a$.\hb
 Then for every $n<\ome$, pick some $\xi_n$ $\del
(p_n)<\xi_n<\kap_n^{+n+2}$  which is $n$-good.
It exists by Lemma 2.3.  Extend $p_{n0}$  to a
condition $q_{n0}$  obtained by adding $\xi_n$
and some $\xi$  which is above $\xi_n$  and
$mc(p_n)$ to $\supp(p_{n0})$.  Set
$q_{n1}=p_{n1}\cup \{\langle \bet ,\xi_n\rangle
\}$,  $q_{n2}=p_{n2}$  and $q_n=\langle q_{n0},q_{n1},
q_{n2}\rangle$.  Then $q=\langle q_n\mid n<\ome\rangle$ 
will be as desired.

\subheading{Case 2}  $\bet <\cup a$.\hb
Then pick the least $\alp\in a$ $\alp >\bet$.
By the definition of $\calP^*$, namely (2) of
2.9, $\alp\in \dom p_{n1}$  starting with some
$n^* <\ome$. by 2.9(3) there is a nondecreasing
converging to infinity sequence of natural
numbers $\langle k_m\mid\ome >m\ge n^*\rangle$
such that for every $m\ge n^*$  $p_{m,1}(\alp)$
is $k_m$-good.  Let $n^{**}\ge n^*$  be such
that $k_{n^{**}}>0$. For every $n\ge n^{**}$ we
like to extend $p_n$  in order to include $\bet$
into the extension.  So, let $n\ge n^{**}$.
Set $\gam=\cup \{p_{n2}(\del)\mid\del<\alp\}$.
Since $p_{n1}(\alp)$ is good.  $cf p_{n1}(\alp)
>\kap_n^{++}$ and hence $\gam<p_{n1}(\alp)$.  by
Lemma 2.5.1, there $k_n-1$-good $\del$,  $\gam
<\del <p_{n1}(\alp)$.  Extend $p_{n0}$  to some
$q_{n0}$  having $\del$ in support.  Set
$q_{n1}=p_{n1}\cup \{ \langle
\bet,\del\rangle\}$,  $q_{n2}=p_{n2}$  and
$q_n=\langle q_{n0}, q_{n1},q_{n2}\rangle$.    

Now for every $n\ge n^{**}$  $q_{n1}(\bet)$ will
be $k_n-1$-good.  Clearly, $\langle k_n-1\mid
n\ge n^{**}\rangle$  is nondecreasing sequence
converging to infinity.  So $q=\langle q_n\mid
n<\ome\rangle$ is a condition in $\calP^*$  as
desired.\hfill$\square$

$\calP^*$  still collapses $\kap_\ome^{++}$  to
$\kap_\ome^+$.  The reason of this as those
of Lemma 1.15.

\proclaim Lemma 2.16.  In $V[G]$
$\ |(\kap_\ome^{++})^\vee |=\kap^+_\ome$.

The following lemma will be the key lemma for
defining the projection of $\calP^*$
satisfying $\kap_\ome^{++}$-c.c. in the next
section.    

But first a definition.

\subheading{Definition 2.17}  Let $p=\langle
p_n\mid n<\ome\rangle$, $q=\langle q_n\mid n <\ome\rangle$
be two conditions in $\calP^*$.  They are called
similar iff
\item{(1)} $\ell (p)=\ell (q)$
\item{(2)} for every $n<\ell (p)$  the following holds  
\item{(2a)} $p_{n0}=q_{n0}$
\item{(2b)} $\min (\dom q_{n1}\bks (\dom
q_{n1}\cap\dom p_{n1}))> \bigcup_{n<\ome}\sup
(\dom p_{n1})$
\item{(2c)} for every $\bet\in \dom p_{n1}\cap
\dom q_{n1}\ p_{n1}(\bet)=q_{n1}(\bet)$ 
\item{(2d)} $|p_{n1}|=|q_{n1}|$  where
$p_n=\langle p_{n0},p_{n1}\rangle$,
$q_n=\langle q_{n0},q_{n1}\rangle$
\item{(3)} for every $n\ge\ell (p)$  the
following holds
\item{(3a)} $p_{n0}=q_{n0}$\hb
 for every $j\in\{ 1,2\}$ 
\item{(3b)} $\min (\dom q_{nj}\bks (\dom
q_{nj}\cap \dom p_{nj}))>\bigcup_{n<\ome}\sup
(\dom p_{nj})$
\item{(3c)} for every $\bet\in\dom p_{nj}\cap
\dom q_{nj}$ $p_{nj}(\bet)=q_{nj}(\bet)$
\item{(3d)} $|p_{nj}|=|q_{nj}|$
where $p_n=\langle p_{n0}, p_{n1},
p_{n2}\rangle$ and $q_n=\langle q_{n0}, q_{n1},
q_{n2}\rangle$.

\proclaim Lemma 2.18.  Suppose $p$  and $q$  are
similar conditions.  Then there are $s\ge p$
and $t\ge q$  such that
\item{(1)} $\ell (s)=\ell (t)$  and $s\rhookup
\ell (s)=t\rhookup \ell (t)$ 
\item{(2)} for every $n\ge \ell(s)$  the
following holds
\item{(2a)} $\dom s_{n1}=\dom t_{n1}=\dom
p_{n1}\cup \dom q_{n1}$
\item{(2b)} $s_{n2}=t_{n2}=p_{n2}\cup q_{n2}$
\item{(2c)} for every $\bet\in\dom s_{n1}=\dom
t_{n1}$  $mc(s_{n0})$ projects to
$s_{n1}(\bet)$ exactly in the same way as
$mc(t_{n0})$ projects to $t_{n1}(\bet)$ 
\item{(3)} there exists a nondecreasing
converging to infinity sequence of natural numbers
$\langle k_n\mid n\ge\ell (s)\rangle$  with
$k_{\ell(s)}\ge 2$  such for every $n\ge \ell
(s)$  the $\calL_{n,k_n,\rho_n}$-type realized
by $mc(s_n)$  and $mc(t_n)$  are identical,
where $\rho_n$  the least upper bound of or the
code of $p_{n1}\tagg (\dom p_{n1}\cap \dom q_{n1})$.

Moreover, if in addition $\min (\bigcup_{\ell
(q)\le n<\ome}\dom q_{n1})\bks \bigcup_{\ell
(q)\le n<\ome}(\dom p_{n1}\cap \dom q_{n1})$  is
in $\dom q_{\ell(q),1}$,  then $s\ge^{*}\!\!p$,
$t\ge^{*}\!\!q$.

\pr Let $\bet$  be the least element of
$\left(\bigcup_{\ell (q)\le n<\ome}\dom q_{n1}\right)
\bks\bigcup_{\ell (q)\le n<\ome}\left(\dom p_{n1}\cap
\dom q_{n1}\right)$.  Pick some $n^*,\ome
>n^*\ge\ell (q)$ such that $\bet\in\dom
q_{n^*,1}$  and for every $n\ge n^*$
$q_{n,1}(\bet)$  is at least $5$-good.
In order to obtain $s$  and $t$  we first extend
$p,q$  to $p',q'$  by adding Prikry sequence up
to level $n^*-1$  such that $\ell (p')=\ell
(q')=n^*$,  $p'\rhookup n^*=q'\rhookup n^*$  and
$p'\bks n^*=p\bks n^*$,  $q'\bks n^*=q\bks n^*$.
Then we apply Lemma 2.8.1.  for every $n$, $\ome
>n\ge n^*$  to $\bet,q'_n$  and $p'_n$ to
produce $t_n$  and $s_n$.  Finally, $t=p'\rhookup 
n^*{}^\cap \langle t_n\mid \ome >n\ge
n^*\rangle$  and $s=p'\rhookup n^*{}^\cap
\langle s_n\mid \ome >n\ge n^*\rangle$  will be
as required.\hfill$\square$

The standard $\Del$-system argument gives the
following

\proclaim Lemma 2.19.  Among any
$\kap^{++}_\ome$-conditions in $\calP^*$  there are
$\kap_\ome^{++}$ which are alike.  

\sect{3.~~The Projection} 

Our aim will be to project $\calP^*$  to a
forcing notion satisfying $\kap^{++}_\ome$-c.c.
but still producing $\kap_\ome^{++}$-Prikry
sequences.

\subheading{Definition 3.0}  Let $n<\ome$  and
suppose $\langle p,f\rangle$, $\langle q,
g\rangle \in Q^*_n$  are such that $f=g$  then
we call them $k$-equivalent for every $k\le n$
and denote this by $\longleftrightarrow_{n,k}$. 

\subheading{Definition 3.1}  Let $2\le k\le n
<\ome$.  Suppose $\langle p,a,f\rangle$, $\langle
q,b,g\rangle \in Q^*_n$.  We call $\langle
p,a,f\rangle$  and $\langle q,b,g\rangle$
$k$-equivalent and denote this by
$\longleftrightarrow_{n,k}$ iff 
\item{(0)} $f=g$
\item {(1)} $\dom a=\dom b$
\item{(2)} $mc(p)$  and $mc(q)$  are realizing
the same $k$-type
\item{(3)} $T(p)=T(q)$, i.e. the sets of measure
1 are the same 
\item{(4)} for every $\del\in \dom a=\dom b\ a(\del)$ 
and $b(\del)$  are realizing the same $k$-type
\item{(5)} for every $\del\in \dom a=\dom b$
and $\ell\le k$  $a(\del)$ is $\ell$-good iff
$b(\del)$  is $\ell$-good 
\item{(6)} for every $\del\in \dom a=\dom b$
$mc(p)$  projects to $a(\del)$  the same way as
$mc(q)$  projects to $b(\del)$.

\subheading{Definition 3.2} Let $p=\langle p_n\mid
n<\ome\rangle$,  $q=\langle q_n\mid
n<\ome\rangle\in\calP^*$.  We call $p$  and $q$
equivalent and denote this by $\longleftrightarrow$
iff
\item{(1)} $\ell (p)=\ell (q)$
\item{(2)} for every $n<\ell (p)$ 
$p_n\longleftrightarrow_{n,n}q_n$, i.e. $p_{n1}=q_{n1}$,
where $p_n=\langle p_{n0},p_{n1}\rangle$  and
$q_n=\langle q_{n0},q_{n1}\rangle$.

Notice that we require only the parts producing
the function from $\kap_\ome^{++}$ to be equal.
So, actually the finite portions of the Prikry
type forcing become unessential. 

\item{(3)} there is a nondecreasing sequence
$\langle k_n\mid \ell (p)\le n<\ome\rangle$,
$\lim_{n\to\infty}k_n=\infty$,  $k_0\ge 2$  such
that for every $n,\ell (p)\le n <\ome$ $p_n$ and
$q_n$  are $k_n$-equivalent.

It is easy to check that $\longleftrightarrow$
is an equivalence relation. 

Now paraphrasing Lemma 2.18 we obtain the
following

\proclaim Lemma 3.3.  Suppose that $p$  and $q$
are similar.  Then there are equivalent $s$
and $t$  such that $s\ge p$  and $t\ge q$.

Note that for every $n\ge\ell (s)=\ell (t)$
$mc(s_{n0})$,  $mc(t_{n0})$ are realizing the
same $\calL_{n,k_n}$-type for $k_n\ge 2$,  where
$s,t$  are produced by Lemma 2.18.  There are at
most $\kap_n^{++}$ different measures over
$\kap_n$.  So, the measures corresponding
$mc(s_{n0})$  and $mc(t_{n0})$ are the same.
Now we can shrink sets of measure one
$T(s_{n0})$  and $T(t_{n0})$ to the same set in
order to satisfy the condition (3) of Definition
3.1.

\subheading{Definition 3.4}  Let $p,q\in\calP^*$.
Then $p\longrightarrow q$ iff there is a
sequence of conditions $\langle r_k\mid
k<m<\ome\rangle$  so that 
\item{(1)} $r_0=p$
\item{(2)} $r_{m-1}=q$
\item{(3)} for every $k<m-1$
$$r_k\le r_{k+1}\quad\hbox{or}\quad r_k\longleftrightarrow
r_{k+1}\ .$$
See diagram:
$$\matrix{r_{m-2}&\longleftrightarrow&r_{m-1}=q\cr
\vee|&&\cr
r_{m-3}&\longleftrightarrow&r_{m-4}\cr
&\cdots&\cr
r_4&\longleftrightarrow&r_5\cr
\vee|&&\cr
r_3&\longleftrightarrow&r_2\cr
&&\vee|\cr
p=r_0&\longleftrightarrow&r_1\cr}$$
Obviously, $\longrightarrow$ is reflexive and
transitive.

\proclaim Lemma 3.5.  Suppose $p,q,s\in\calP^*$
$p\longleftrightarrow q$  and $s\ge p$.  Then
there are $s'\ge s$  and $t\ge q$  such that
$s'\longleftrightarrow t$.

\pr Pick a nondecreasing sequence $\langle
k_n\mid \ell (p)=\ell(q)\le n<\ome\rangle$,
$\lim_{n\to\infty}k_n=\infty$  such that
$p_n\longleftrightarrow_{n,k_n}q_n$  for every
$n\ge\ell(p)$.  For each $n$,  $\ell (p)\le n
<\ell (s)$ we extend $q_n=\langle q_{n0},q_{n1},q_{n2}
\rangle$ to $t_n=\langle t_{n0},t_{n1}\rangle$
by putting $s_{n0}^{mc(p_{n0})}$  over $mc(q_{n0})$
projecting it over the rest of the coordinates in $\supp
q_{n0}$  and $rng q_{n1}$  and setting
$t_{n1}=s_{n1}$, where $s_n=\langle
s_{n0},s_{n1}\rangle$,  $p_n=\langle
p_{n0},p_{n1},p_{n2}\rangle$  and
$s_{n0}^{mc(p_{n0})}$ is the one element
sequence standing over the maximal coordinate of
$p_{n0}$.  Notice that this is possible since $T(p_{n0})=
T(q_{n0})$  and $s_{n0}^{m(p_{n0})}\in T(p_{n0})$.  Then
$s_n$  and $t_n$  will be $n$-equivalent.  Set
$s'_n=s_n$.

Suppose now that $n\ge \ell (s)$.  Let
$s_n=\langle s_{n0}, s_{n1},s_{n2}\rangle$,
$p_n=\langle p_{n0},p_{n1}p_{n2}\rangle$  and
$q_n=\langle q_{n0},q_{n1},q_{n2}\rangle$. 

\subheading{Case 1} $k_n>2$.\hb
By Lemma 2.0, there is $\del$  realizing the same
$k_n-1$-type over $mc(q_{n0})$  as $mc(s_{n0})$
does over $mc(p_{n0})$.  Now pick  $t_n=\langle
t_{n0},t_{n1},t_{n2}\rangle$  to be a condition
with $mc(t_{n0})=\del$  $k_n-1$-equivalent to
$s_n$. Set $s'_n=s_n$.

\subheading{Case 2} $k_n\le 2$.\hb  
We first extend $s_n$  to a stronger condition
$s'_n=\langle s'_{n0},s'_{n1}\rangle$.  Then we
proceed as in the case $\ell (p)\le n <\ell
(s)$.

By the construction $s'=\langle s'_n\mid
n<\ome\rangle$ and $t=\langle t_n\mid
n<\ome\rangle$  will be stronger than $s$ and
$q$  respectively.  Also $\ell(s')=\ell (t)$ and
for every $n<\ell (s)$  $s'_n\longleftrightarrow_{n,n}$
$t_n$.  The sequence $\langle k_n-1\mid \ell
(s')\le n<\ome\rangle$  will witness the
condition (2) of Definition 3.2.\hfill$\square$

Now let us define the projection.

\subheading{Definition 3.5} Set 
$$\calP^{**}=\calP/\longleftrightarrow\ .$$

For $x,y\in\calP^{**}$  let $x\preceq y$  iff
there are $p\in x$  and $q\in y$  such that 
$p\longrightarrow q$.

\proclaim Lemma 3.7.  A function $\pi:\calP^*\to
\calP^{**}$ defined by $\pi (p)=p/\longleftrightarrow$
projects $\langle \calP^*,\le\rangle$  nicely
onto $\langle \calP^{**},\preceq \rangle$.

\pr It is enough to show that for every $p,q\in
\calP^*$ if $p\to q$ then there is
$s\ge p$  such that $q\to s$.  Suppose for
simplicity that we have the following diagram
witnessing $p\to q$.  In a general case the same
argument should be applied inductively.
$$\matrix{q&\longleftrightarrow&h\cr
&&\vee|\cr
f&\longleftrightarrow&g\cr
\vee|&&\cr
d&\longleftrightarrow&c\cr
&&\vee|\cr
a&\longleftrightarrow&b\cr
\vee|&&\cr
p&&\cr}$$
Using Lemma 3.5 we find equivalent $f'\ge f$
and $h'\ge h$.  Then applying it to $d,c,f'$
find equivalent $f\tagg \ge f'$ and $c\tagg\ge c$. 
Finally, using Lemma 3.5 for $c\tagg$,  $b,a$
we find equivalent $a\taggg \ge a$  and
$c\taggg\ge c\tagg$.  In the diagram it looks
like:
$$\matrix{&&&&&&h'\cr
&&&&&&\vee|\cr
&&&&q&\longleftrightarrow&h&&&&\cr
&&&&&&\vee|\cr
f\tagg&\ge&f'&\ge&f&\longleftrightarrow&g&&&&\cr
&&&&\vee|&&&&&&\cr
&&&&d&\longleftrightarrow&c&\le&c\tagg&\le&c\taggg\cr
&&&&&&\vee|&&&&\cr
&&a\taggg&\ge&a&\longleftrightarrow&b&&&&\cr
&&&&\vee|&&&&&&\cr
&&&&p&&&&&&\cr}$$

We claim that $a\taggg$  is as required, i.e.
$a\taggg\ge p$ and $q\longrightarrow a\taggg$.  Clearly,
$a\taggg\ge p$. In order to prove $q\longrightarrow
a\taggg$  we consider the following diagram:
$$\matrix{a\taggg&\longleftrightarrow&c\taggg\cr
&&\vee|\cr
f\tagg&\longleftrightarrow&c\tagg\cr
\vee|&&\cr
f'&\longleftrightarrow&h'\cr 
&&\vee|\cr
q&\longleftrightarrow&h\cr}$$
So the sequence $\langle
q,h,h',f',f\tagg,c\tagg,c\taggg,a\taggg\rangle$
witnessing $q\longrightarrow a\taggg$.\hfill$\square$

The next lemma follows from Lemma 3.3. 

\proclaim Lemma 3.8.  $\calP^{**}$  satisfies
$\kap_\ome^{++}$-c.c.

Let $G\subseteq\calP^*$  be generic.  We like to
show that for every $\bet <\kap_\ome^{++}$
$G(\bet)\in V[\pi\tagg (G)]$.

The following will be sufficient.

\proclaim Lemma 3.9. Let $p\longleftrightarrow
q$,  $\bet <\kap_\ome^{++}$. Suppose that for
some $n<\ell(p)$  $\bet\in\dom p_{n1}$  then
$\bet \in \dom q_{n1}$  and
$p_{n1}(\bet)=q_{n1}(\bet)$.  Where $p_n=\langle
p_{n0},p_{n1}\rangle$  and $q_n=\langle q_{n0},q_{n1}
\rangle$.

\pr By the definition of equivalence
$q_{n1}=p_{n1}$.\hfill$\square$

So using Lemma 2.14 we obtain the following

\proclaim Theorem 3.10.  Let $G$  be a generic
subset of $\calP^*$.  Then $V[\pi\tagg (G)]$ is
a cardinal preserving extension of $V$  such
that GCH holds below $\kap_\ome$  and
$2^{\kap_\ome}=\kap_\ome^{++}$. 

\sect{4.~~Down to $\aleph_\ome$}

In this section we sketch an additional
construction needed for moving $\kap_\ome$  to
$\aleph_\ome$.  The construction will be similar
to those of [Git-Mag1]. 

Let $G$  be a generic subset of the forcing
$\calP^{**}$  of the previous section.  Denote
by $\langle\rho_n\mid n<\ome \rangle$  a Prikry
sequence corresponding to normal measures over
$\kap_n$'s.  Then
$cf\left(\prod_{n<\ome}\rho_n^{+n+2}/\hbox{finite}\right)
=\kap_\ome^{++}$. 
Just $G(\bet)$'s $(\bet <\kap_\ome^{++})$  which
are Prikry sequences are witnessing this.
The idea will be to collapse $\rho_{n+1}$ to
$\kap_n^{+n+2}$  and all the cardinals between
$\rho_{n+1}^{+n+4}$  and $\kap_{n+1}$  to
$\rho_{n+1}^{+n+4}$.  In order to perform this
avoiding collapse of $\kap_\ome^{++}$,  we need
modify $\calP^*$.  For collapsing cardinals
between $\rho_{n+1}^{+n+4}$  and $\kap_{n+1}$
the method used in [Git-Mag 1] applies directly
since the length of the extender used over $\kap_{n+1}$
is only $\kap_{n+1}^{+(n+1)+2}$.  Hence let us
describe only the way $\rho_{n+1}$ will be collapsed
to $\kap_n^{+n+2}$.

Let us deal with a fixed $n<\ome$  and drop the
lower index $n$ for a while.  Fix a nonstationary
set $A\subseteq \kap^{+n+2}$.  In Definition 1.2
we require in addition that $rng\cap
A=\emptyset$  and $\supp p\cap A=\emptyset$.  In
the definition of the order on $Q$, Definition
1.4 (2) for $\gam\in A$ we replace $p^\gam$  by
$\kap$  only if $p^\gam\ge\kap_{n+1}$.  Now, the
definition of $\calP$, Definition 1.8 is changed
as follows:

\subheading{Definition 4.1}  A set of forcing
conditions $\calP$  consists of all elements $p$
of the form $\langle p_n\mid n<\ome\rangle$  so
that 
\item{(1)} for every $n<\ome$ $p_n\in Q_n$ 
\item{(2)} there exists $\ell <\ome$  such that
for every $n\ge \ell$  $p_n\in Q^0_n$
\item{(3)} if $0<n<\ell (p)$,  then for every
$\gam\in A_{n-1}\cap \del (p_{n-1,0})$  $p^\gam_{n-1,0}<
p^0_{n,0}$, where $p_n=\langle
p_{n0},p_{n1}\rangle$ and  $p_{n-1}=\langle p_{n-1,0},
p_{n-1,1}\rangle$. 

The meaning of the new condition (3) is that
$p^0_{n0}$  which is $\rho_n$  is always above
all the sequences mentioned in $p_{n-1,0}$.
This will actually produce a cofinal function
from $A_n$  into $\rho_n$.

Finally, in order to keep it while going to the
projection $\calP^{**}$, we strengthen the notion
of similarity.  Thus, in Definition 2.17 we
require in addition that for every $\gam\in
a_n\cap \del (p_{n0})\ p^\gam_{n0}=q^\gam_{n0}$.
I.e.  the values of the cofinal function $A_n\mapsto
\rho_n$  are never changed.

There is no problem in showing the Prikry
condition, (i.e. Lemma 1.11) since passing from
level $n-1$ to level $n$ we will have a regressive
function on a set of
measure one for a normal measure over $\kap_n$. 

\sect{5.~~Loose Ends}

We do not know if it is possible under the same
initial assumption to make a gap between
$\kap_\ome$  and $2^{\kap_\ome}$  wider.
Our conjecture is that it is possible.  Namely,
it is possible to obtain countable gaps.  Also
we think that uncountable gaps are
impossible.
\vskip1truecm
\references {60}
\ref{[Git]} M. Gitik, On Hidden Extenders
\smallskip
\ref{[Git-Mit]} M. Gitik and W. Mitchell,
Indiscernible Sequences for Extenders and the
Singular Cardinal Hypothesis.
\smallskip
\ref{[Git-Mag1]} M. Gitik and M. Magidor, The
Singular Cardinal Hypothesis Revisited, in MSRI
Conf. Proc., 1991, 243-279.
\smallskip
\ref{[Git-Mag2]} M. Gitik and M. Magidor,
Extender Based Forcing Notions, to appear in
JSL. 
\smallskip
\ref{[Mit-St-Sch]} W. Mitchell, J. Steel and E.
Schimmerling.

\end

%% file: def1200.tex
%


\def\today{\ifcase\month\or January\or February\or
March\or April\or May\or June\or July\or August\or
September\or October\or November\or December\fi
\space\number\day, \number\year}




\def\dspace{\lineskip=2pt\baselineskip=18pt
\lineskiplimit=0pt}

\font \bbrm=cmbx10 at 12pt

\def\bigtype{\bbrm}

\hsize=13.5cm
\magnification=1200
\def\ce{\centerline}

\def\hb{\hfill\break}

\def\title #1{\null\bigskip\ce{\bigtype #1}
\bigskip}

\def\alp{\alpha}		
\def\bet{\beta}		
\def\gam{\gamma}		
\def\del{\delta}		\def\Del{\Delta}

\def\kap{\kappa}
\def\lam{\lambda}		
\def\sig{\sigma}		

\def\ome{\omega}		


\def\calL{{\cal L}}

\def\calP{{\cal P}}

\def\calU{{\cal U}}



    
\font\tenboldgreek=cmmib10
 \font\sevenboldgreek=cmmib10 at 7pt
\font\fiveboldgreek=cmmib10 at 7pt
\newfam\bgfam
\textfont\bgfam=\tenboldgreek
\scriptfont\bgfam=\sevenboldgreek
\scriptscriptfont\bgfam=\fiveboldgreek

\mathchardef\ggarrow="7010

\font\tengerman=eufm10 \font\sevengerman=eufm7
\font\fivegerman=eufm5
\font\tendouble=msym10 \font\sevendouble=msym7
\font\fivedouble=msym5

\textfont4=\tengerman \scriptfont4=\sevengerman
\scriptscriptfont4=\fivegerman
\newfam\dbfam
\textfont\dbfam=\tendouble \scriptfont\dbfam=
\sevendouble
\scriptscriptfont\dbfam=\fivedouble
\def\gr{\fam4}

\mathchardef\ng="702D
\mathchardef\dbA="7041
\mathchardef\sm="7072
\mathchardef\nvdash="7030
\mathchardef\nldash="7031
\mathchardef\lne="7008
\mathchardef\sneq="7024
\mathchardef\spneq="7025
\mathchardef\sne="7028
\mathchardef\spne="7029
\mathchardef\ltms="706E
\mathchardef\tmsl="706F

\mathchardef\dbA="7041

	\def\gra{{\gr a}}

\mathchardef\dbA="7041 
\mathchardef\dbB="7042 
\mathchardef\dbC="7043 
\mathchardef\dbD="7044 
\mathchardef\dbE="7045 
\mathchardef\dbF="7046 
\mathchardef\dbG="7047 
\mathchardef\dbH="7048 
\mathchardef\dbI="7049 
\mathchardef\dbJ="704A 
\mathchardef\dbK="704B 
\mathchardef\dbL="704C 
\mathchardef\dbM="704D 
\mathchardef\dbN="704E 
\mathchardef\dbO="704F 
\mathchardef\dbP="7050 
\mathchardef\dbQ="7051 
\mathchardef\dbR="7052 
\mathchardef\dbS="7053 
\mathchardef\dbT="7054 
\mathchardef\dbU="7055 \def\UU{{\fam=\dbfam\dbU}}
\mathchardef\dbV="7056 
\mathchardef\dbW="7057 
\mathchardef\dbX="7058 
\mathchardef\dbY="7059 
\mathchardef\dbZ="705A 

\def\nek{,\ldots,}
\def\sdp{\times \hskip -0.3em {\raise 0.3ex
\hbox{$\scriptscriptstyle |$}}} 


\def\dom{\mathop{\rm dom}\nolimits}

\def\min{\mathop{\rm min}}
\def\MOD{\mathop{\rm mod}}

\def\supp{\mathop{\rm supp}}









\def\tileta{{\widetilde\eta}}

\def\ddownarrow{\big\downarrow \hskip-0.70em\raise
2pt\hbox {$\big\downarrow$}}
\def\longright #1#2 {\smash{\mathop{\hbox to
#1pt {\rightarrowfill}}\limits_{#2}}}
\def\sqr#1#2{{\vcenter{\hrule height.#2pt\hbox{\vrule
width.#2pt height#1pt \kern#1pt \vrule width.#2pt}
\hrule height.#2pt}}}
\def\square{\mathchoice{\sqr34}{\sqr34}{\sqr{2.1}3}
{\sqr{1.5}3}}

\def\buildrul#1\under#2{\mathrel{\mathop{\null#2}
\limits_{#1}}}

\def\boxit#1{\vbox{\hrule\hbox{\vrule\kern3pt
\vbox{\kern3pt#1 \kern3pt}\kern3pt\vrule}\hrule}}

\def\bigcupl{\bigcup\limits}

\def\subheading#1{\medskip\goodbreak\noindent{\bf
#1.}\quad}

\def\sect#1{\goodbreak\bigskip\centerline{\bf#1}
\medskip}
\def\pr{\smallskip\noindent{\bf Proof:\quad}}
\def\onumber #1{\ooalign{\hfil\raise.07ex\hbox{
\hfill$\scriptstyle \,#1$\hfil}
\cr\cr{$\bigcirc$}}}
\def\onumber c{\ooalign{\hfil\raise.07ex\hbox
{\hfill$\scriptstyle \,c$\hfil}
\cr\cr{$\bigcirc$}}}
\def\alpcirc {\ooalign{\hfil\raise.07ex
\hbox{\hfill$\scriptstyle\alp\;$\hfill}\cr\cr
{$\bigcirc$}}}

\def\longmapright #1#2 {\smash{\mathop{\hbox to
#1pt {\rightarrowfill}}\limits^{#2}}}
\def\longmapleft #1 #2 {\smash{\mathop{\hbox to
#1 pt {\leftarrowfill}}\limits^{#2}}}

\def\references#1{\goodbreak\bigskip\par\centerline
{\bf References}\medskip\parindent=#1pt}
\def\ref#1{\par\smallskip\hang\indent\llap{\hbox
to \parindent{#1\hfil\enspace}}\ignorespaces}

\def\back{{\raise 2.5pt\hbox{$\,\scriptscriptstyle
\backslash\,$}}}
\def\bks{{\backslash}}
\def\part{\partial}
\def\lwr #1{\lower 5pt\hbox{$#1$}\hskip -3pt}
\def\rse #1{\hskip -3pt\raise 5pt\hbox{$#1$}}
\def\lwrs #1{\lower 4pt\hbox{$\scriptstyle #1$}
\hskip -2pt}
\def\rses #1{\hskip -2pt\raise 3pt\hbox
{$\scriptstyle #1$}}

\def\<#1{\left\langle{#1}\right\rangle}

\def\subinbn{{\subset\hskip-8pt\raise 0.95pt
\hbox{$\scriptscriptstyle\subset$}}}

\def\llvdash{\mathop{\|\hskip-2pt
\raise 3pt\hbox{\vrule height 0.25pt width 1.5cm}}}

\def\lvdash{\mathop{|\hskip-2pt \raise 3pt\hbox
{\vrule height 0.25pt width 1.5cm}}}

\def\fakebold#1{\leavevmode\setbox0=\hbox{#1}%
  \kern-.025em\copy0 \kern-\wd0
  \kern .025em\copy0 \kern-\wd0
  \kern-.025em\raise.0333em\box0 }

\font\msxmten=msxm10
\font\msxmseven=msxm7
\font\msxmfive=msxm5
\newfam\myfam
\textfont\myfam=\msxmten
\scriptfont\myfam=\msxmseven
\scriptscriptfont\myfam=\msxmfive
\mathchardef\rhookupone="7016
\mathchardef\ldh="700D
\mathchardef\leg="7053
\mathchardef\ANG="705E
\mathchardef\lcu="7070
\mathchardef\rcu="7071
\mathchardef\leseq="7035
\mathchardef\qeeg="703D
\mathchardef\qeel="7036
\mathchardef\blackbox="7004
\mathchardef\bbx="7003
\mathchardef\simsucc="7025

\def\rhookup{{\fam=\myfam \rhookupone}}

\font\tencaps=cmcsc10
\def\smallcaps{\tencaps}

\def\author#1{\bigskip\ce{\smallcaps #1}\medskip}

\def\tagg{^{\prime\prime}}
\def\taggg{^{\prime\prime\prime}}

\def\upddots{\mathinner{\mkern
1mu\raise 1pt \hbox{.}\mkern 2mu \mkern
2mu \raise 4pt\hbox{.}\mkern 1mu \raise 7pt\vbox
{\kern 7 pt\hbox{.}}} }

\def\varchi{\ooalign{{\raise
1.385pt\hbox{$\chi$}}\crcr\hbox{--}\crcr}}

\def\trianarrow{{\raise 2pt\hbox to 0.50cm
{\hrulefill}\triangleright}}